\documentclass[12pt]{article}
\usepackage{geometry}
\usepackage{amssymb}
\usepackage{amsmath}
\usepackage{graphicx}
\usepackage{endnotes}
\usepackage{array}
\usepackage{url}
\usepackage{sectsty}

\usepackage[backref=page, 
hyperfootnotes=true 
]{hyperref} 

\hypersetup{
colorlinks=true, 
linkcolor=red, 
anchorcolor=black, 
citecolor=blue, 
urlcolor=cyan, 
bookmarksnumbered=true, 
linktocpage=false 
}

\newcommand\vel{\vee}
\newcommand\ts{\vdash}

\newcommand{\onlyif}{\rightarrow}
\newcommand{\Russell}{\rotatebox[origin=c]{180}{$\iota$}}  
\title{\textbf{Identity, Haecceity, and the Godzilla Problem}\footnote{Published in Gillman Payette (ed.), \textsl{``Shut Up'' he explained.  Essays in Honour of Peter K. Schotch}.  London:  College Publishers, 2016, pp.\ 63--79.}}
\date{}
\author{Kent A. Peacock \\ \small{\textsl{Department of Philosophy, University of Lethbridge}} \\ \small{\url{kent.peacock@uleth.ca}}
\and Andrew Tedder \\ \small{\textsl{Department of Philosophy, University of Connecticut}} \\ \small{\url{andrew.tedder@uconn.edu}}
}  

\begin{document}
\maketitle
\linespread{1.25}\selectfont\raggedright\parindent20pt
\allsectionsfont{\normalsize}

\begin{abstract}
  \noindent In standard first order predicate logic with identity it is
  usually taken that $a=a$ is a theorem for any term $a$.  It is
  easily shown that this enables the apparent proof of a theorem stating the existence of any
  entity whatsoever.  This embarrassing result is a motivation for the construction of free
  logics, but in most orthodox treatments of first order logic with identity it is generally dealt
  with by being ignored.  We investigate the possibility that this
  problem can be obviated by dropping the rule that $a=a$ is a
  theorem and requiring instead that it be treated as a global but in principle defeasible
  assumption about the objects in the domain of discourse.  
  We   review some motivations in physics, philosophy, and
  literature for questioning the classical notion of self-identity, and we show that Carnap's ``null object'' has a natural role to play in any system of predicate logic where self-identity can come into question.
\end{abstract}

\section{Introduction:  An Embarrassing Problem}
If one's only knowledge of logic came from standard university texts, one
might think that elementary first order predicate logic with identity has all
been worked out a long time ago and that there are thus no serious technical
or conceptual problems lurking within it. We're not sure that this
comfortable view is right.

We're going to begin by pointing to what we believe is an embarrassing
problem for standard first-order predicate logic with identity. The usual
approach is to take the self-identity of all objects in the universe of
discourse as a logical truth or theorem; that is, it is taken that
\begin{equation}
\ts \forall x(x=x); \label{IdTh}
\end{equation}
or as
\begin{equation}
  \ts a=a
\end{equation}
in a scheme allowing generalization to (\ref{IdTh}). The symbol $\ts$ in
these formulas is not meant to suggest that they are provable from the other
resources of first-order logic; but rather that they are assertions that may be 
added without proof to
first order predicate logic to give a theory of identity.  In Lemmon's system
\cite{Lemmon65} the rule (\ref{IdTh}) is called Identity Introduction and symbolized
=I. It is paired with Identity Elimination (=E), which says that if $a=b$
then $a$ may be substituted for $b$ or \textit{vice versa} wherever they
occur. Identity Elimination is the motor that drives virtually
all of all applications of identity.  Identity Intro, however, is rarely
used, and as the following natural deduction proof shows, it allows
consequences with which all logicians ought to be uncomfortable:
\begin{tabbing}
mm \= mmmm \= mmm  \= mmmmmmmmmmmmmm  \= mmmmmm   \kill
    \>     \>(1)  \> $\forall x(x=x)$           \>=I  \\
    \>     \>(2)  \> $\text{Godzilla} = \text{Godzilla}$      \>1 UE  \\
    \>     \>(3)  \> $\exists x(x = \text{Godzilla})$    \>2 EI
\end{tabbing}
So we have
\begin{equation}
  \ts \exists x(x = \text{Godzilla}).
\end{equation}
We'll call this argument pattern the \emph{Categorical Godzilla}.
(There is also a \emph{Conditional Godzilla} which we will later
introduce.)  If $\forall x(x=x)$ is a logical truth, then we apparently have it as a \emph{logical truth} that a
certain city-trampling movie monster really does exist.  This pattern could be repeated for any name whatsoever.  
\emph{What is wrong with this picture?}

Of course, this problem is well known (or should be) and it is one of the
motivations for the construction of various types of free logic, the
defining characteristic of which is that names are not automatically
assumed to refer \cite{Lambert91,BDS,Nolt11}. We are sympathetic to free
logic and we think that our observations here tend to give further
motivation for its development.  However, our primary aim in this
note is the more limited goal of arguing against the theoremhood of
=I.  In doing so, we will show that there is reason to question =I even as a universal assumption for at least some of the domains to which predicate logic might be applied.
We will also show that a consequence of treating =I as a global but defeasible assumption is, surprisngly perhaps, the resurrection of Carnap's \emph{null thing}  \cite{Carnap47}, or something very much like it.

\section{A Short History of Arguments for =I}
Many currently used logic texts get around the awkward conclusion that =I can
be used to prove that anything whatsoever exists simply by not mentioning it
at all. E.\ J.\ Lemmon's widely used \textsl{Beginning Logic} \cite{Lemmon65}
relegates it without comment to an exercise, and defends the theoremhood of
=I as follows:
\begin{quote}
  For any term $t$, the rule =I permits us to introduce into a proof
  at any stage $t=t$, resting on no assumptions.  The idea should be
  clear:  anything is itself, as a matter of logic; hence $t=t$ is
  logically true, and so can appear without assumptions. \cite[p.\
  161]{Lemmon65}
\end{quote}
It is by no means clear that everything being itself is a matter of logic. By
contrast, the theorem $P \onlyif P$ is a matter of logic, and indeed it is
sometimes (confusingly) called Identity; it is almost as if Lemmon confused
$=$ with the truth-functional connective $\onlyif$, or perhaps $\equiv$.

A recent logic text by Paul Herrick trades on intuitions similar to
Lemmon's:
\begin{quote}
  Each thing is identical with itself \dots Surely this needs no
  argument; certainly it is necessarily true.  (How could something
  possibly \emph{not} be identical to itself?) \cite[p.\ 587]{Herrick13}
\end{quote}
One should be suspicious of arguments for $p$ of the form, ``Surely
$p$\dots''.  For surely (if we may) Lemmon and Herrick are appealing
not to logical intuitions about identity, but \emph{metaphysical}
intuitions about identity. They are saying that it is a matter of
\emph{necessary fact} that every item is self-identical, and they
have forgotten that pure logic as such expresses no facts. Facts,
whether necessary or of the ordinary garden-variety, can be
introduced into a logical problem only \emph{by assumption}.  The
views of Lemmon and Herrick therefore seem to be part of a long
tradition of mistaking presumed factual or metaphysical necessity
for logical or mathematical necessity.

In \textit{Principia Mathematica} Russell and Whitehead \cite[$*$13]{PM} gave
an apparently much more principled and precise defence of =I. They begin by
defining $=$ by means of the Identity of Indiscernibles:  $x=y$ means that if
$\phi$ is a property of $x$ then $\phi$ is a property of $y$.  (In this
sketch of their exposition we gloss over niceties  having to do with the
Theory of Types.)  The definition gets turned into a theorem by an
application of universal instantiation:  since it is true of any two
arbitrarily selected entities that they are identical if and only if they
share all properties, then that holds for all instances of $x=y$. Then as a
special case of this result any arbitrarily selected entity is self-identical
simply because any property of itself is a property of itself.

The Russell-Whitehead approach has the virtue of precision.  However, modern
authors tend to shy away from it because it requires second order logic. More
important from the skeptical point of view we pursue here, the reliance on
the Identity of Indiscernibles again amounts to building a metaphysical
principle into formal logic.  It could even be said (though no doubt
contentiously) that Russell and Whitehead's view, that all terms are self-identical because any arbitrary term would have the same properties as itself,
borders on question-begging. Russell and Whitehead do not state the
Categorical Godzilla, though it would be readily available in their system.

Reliance on the Identity of Indiscernibles in order to define identity and
justify theorems about it can be traced back through Frege \cite{Begriff} and
Leibniz (discussed in Kneale \cite{Kneale62}) to Aristotle.  The latter seems
to be the first to have introduced the concept of the Identity of
Indiscernibles, though informally, in \emph{De Sophisticis Elenchis}
\cite[Ch.\ 24 (179a37)]{Aristotle_SE}: ``For only to things that are indistinguishable
and one in essence is it generally agreed that all the same attributes
belong.''  His wording ``it is generally agreed'' suggests that this
principle has a longer history, either orally or in writings now lost.  Being ``one in essence'' is a metaphysical requirement for self-identity; Aristotle's discussion 
surrounding the line quoted here explains why identity may otherwise be ambiguous if this high metaphysical standard is not met.    

The reliance upon the Identity of Indiscernibles to get =I is
therefore very old. Now =I, for our purposes, can be
treated as either substantive or simply stipulative. If the latter,
then it's the kind of thing that we may choose to do without. If the
former, then its place in the reasoning designed to justify logic is
dubious at best. Logic should not be a substantive inquiry, but
rather a formal one (this is how the field has been moving, and we
think, for the better). What we should be interested in are
extremely generalized relations between assertions about objects burdened with as few assumptions as possible.  So how can we amend standard first
order logic with identity in as conservative a way as possible, but
so as to block the Categorical Godzilla?

\section{Blocking the Godzilla Inference}
There is no way to block the Godzilla proof by placing some sort of
artificial restriction on EI without crippling or drastically
reconstruing EI, which is not in line with our conservative
approach. And the application of UE to line (1) seems to be entirely
unobjectionable.  
One could simply introduce an \textit{ad hoc} rule
against making inferences of the Godzilla type; for instance, one
might introduce a rule against applying EI to any formula of the
form $a=a$. This probably would be logically possible but it seems
inelegant. 
To begin with, we want to change as little as possible of classical
first order logic, and only make such changes to identity theory as
would be sufficient to block Godzilla in a natural way.  (Further along we'll suggest something more radical.)  It is
therefore much more pertinent to examine whether it really is
reasonable to treat $\forall x(x=x)$ as a theorem.

The first point to note is that the Godzilla proof is simply an
illustration of the rule that there is an existential claim built
into any proposition of the form $Fa$. Suppose we assume or are
given $Fa$, where $F$ is some property and $a$ is an individual.
From this we conclude $\exists xFx$ by an immediate application of
EI. So to assert that an individual does in fact possess a property
is to imply the existence of an individual possessing that property,
and this is true for any property $F$, \emph{including
self-identity}.

This suggests that another way of blocking the Categorical Godzilla might be to question whether self-identity really can be treated as if it were monadic property.  If Godzilla is a city-trampler then something is a city-trampler; if Godzilla is self-identical then can we say that something is self-identical?  It would be very odd if we could not, since self-identity certainly \emph{is} something that pertains to individual entities when it pertains at all; hence, that approach does not seem promising either.  

The key is that we do not want to build any assumptions about the existence of
any entities whatsoever into our logic. Therefore, we can avoid the Categorical Godzilla if (i) we remind ourselves that
$\forall xFx$ and $Fa$ can be introduced in a proof only as
\emph{assumptions} and (ii) insist that this rule be followed even
when $F$ is self-identity---and even if we find it hard to
\emph{imagine} that any given entity could not be self-identical. On
some metaphysical views it might indeed be the case that everything
is self-identical, but this can't be a matter of \emph{logic} even
if it is a necessary truth, whatever \emph{that} might be.  There
are certainly some universes of discourse that contain only
self-identical objects (such as the universe comprised of the set of natural numbers), but our
choice of a universe of discourse is not a matter of logic either.

A possible defence of the orthodox view could be along the following
lines: it could be said that in doing first-order logic we always
take it for granted that the universe of discourse $\mathcal{U}$
consists of objects that are already presumed to exist.  We see two
objections to this view.

First, even if we want to say that we take it from the outset that
all items in $\mathcal{U}$ exist, we do not mean to say that it
turns out to be a \emph{theorem} that some item in $\mathcal{U}$,
which might be, say, the city of Paris, France, happens to exist;
rather, the existence of any item in the actual world is an
empirical matter if we are talking about real-world entities (unless
we could view the world from the perspective of the God of Leibniz
for whom all apparently empirical truths are analytic), or a
mathematical matter if we are talking about mathematical entities
such as sets or numbers.

Second, logic would not be very useful if we only allowed ourselves
to talk about things that we already know or believe to exist.  One
of the most powerful tools of thought is the ability to consider
hypothetical objects without existential commitment.  This is just
what we do in indirect proof in mathematics. Consider Euclid's proof
that there is no greatest prime:  it begins by supposing
hypothetically that there is a greatest prime, it gives this hypothetical number a name for convenience in calculation, and then shows that
such a number must have contradictory properties. The universe of
discourse should be precisely that---the collection of entities and
subjects that we want to \emph{talk} about and \emph{investigate}, without
necessarily having made a prior commitment to their existence.  We
only get existence out of a deduction if we have good reason to put
it in as a premise.  Of course, we may wish to do logic over
specific sets of entities, such as the natural numbers or
the items of furniture in someone's office, which are already taken
to exist and to have various definite properties.  But logic should
preserve a general freedom to talk of the hypothetical or the
fictional, and to deny the existence of entities if the facts demand that we do so.  

The description of a hypothetical or fictional object may or may not include or imply self-identity.   Logic does not rule out either of these possibilities.  For instance, the notion of the largest prime implies self-identity by the Peano axioms, even though there is no such number.  However, it is not at all clear whether the script writers who defined Godzilla intended to imply that the creature should be self-identical; that doesn't seem to have been relevant.  We are entirely at liberty to define anything with any putative properties whatsover, be it Quine's round square cupola or a non-self identical movie monster.  As Gaunilo realized a long time ago, no definition by itself can have any bearing on what exists.  A definition can neither \emph{make} something exist in the real world (as supposed by the confident authors of various ontological proofs) nor can it, by itself, \emph{prevent} something from existing (even if the description of the putative entity is contradictory).  In the end, as Hume indicated, the test of existence outside the lands of mathematics is always empirical.  

We suggest that we can avoid the Godzilla inference, as least so far
as first order logic is concerned, by taking the following two
steps:
\begin{itemize}
\item We entirely drop the idea that (\ref{IdTh}) is a theorem that applies
    to any arbitrarily selected class of entities. Instead, any
    suggestion that some entity or set of entities are self-identical
    must be treated as an  \emph{assumption} and introduced to a natural
    deduction proof accordingly, unless it is specified in advance that
    one is quantifying over sets of entities (such as the natural
    numbers) which are already known to be self-identical.  In such a
    case, a line in a deduction that says (say) $3=3$ is justified by
    reference to number theory, not to =I.
\item  Within the object language of first order logic we treat identity
    in a purely syntactic manner:  if $a=b$ all that this means is that
    $a$ can be freely substituted for $b$ or \textit{vice versa} at any point in a proof.
\end{itemize}
The key new idea is that we do not take it as a theorem that $\forall
x(x=x)$.  Rather, we take it that identity has to be either postulated
or established from other postulates for any entity or class of entities.  

The purely syntactic reading of $a=b$ reflects the fact that
inter-substitution is all that matters from the syntactic point of view. Most
important, we can define identity this way without having to worry about deep
metaphysical questions about what it means for entities to be identical. This
is hardly to suggest that one should not investigate the nature of identity
or the various sorts of identity that might be tenable (for indeed, our
minimal syntax could be compatible with several interpretations of identity);
it is simply to insist that views about the factual or metaphysical nature of
identity should not be surreptitiously built into first-order logic.  So
self-identity cannot be taken as a given for all possible objects.  Instead,
a self-identity claim is to be introduced into a proof if needed by an
assumption (which could be either universally quantified or of the form $a=a$
for some particular term $a$) or by specifying in advance that one is
quantifying over a class of objects (such as natural numbers) that are
already known to be self-identical.

If we want a system that works the same way as standard first-order predicate logic with identity, so that we can do all of the usual definite description problems and employ the other non-problematic applications of identity, we can take it as not a logical truth but a \emph{global assumption} (which need not be stated explicitly in each deduction) that we reason over domains of objects that are presumed to be self-identical.  The logical status of =I then becomes something like the logical status of the parallel postulate in geometry.  The parallel postulate was once presumed to be either \textit{a priori} or deducible from the other rules of Euclidean geometry, but by the early 19th century it was evident that it was a logically independent assumption about the kinds of spaces that one was dealing with.  (From the Riemannian viewpoint  it applies to spaces with zero intrinsic curvature.)  Logics in which self-identity fails or could fail for some or all non-null entities in the domain of discourse can be called \emph{non-Aristotelian}, by analogy with non-Euclidean geometry.  A logic that differs from classical (Aristotelian) predicate logic with identity \emph{only} in that =I is taken to be a global assumption rather than a logical truth we will call an \emph{open classical logic}---open because it is open to the possibility that the self-identity assumption might fail for one or more members of the domain.  We will show below that even open classical logics must be non-Aristotelian in one particular respect, but there could well be many possible non-Aristotelian logics, just as there are many possible non-Euclidean geometries.  In this respect, logic still seeks its Riemannian synthesis.

\section{Haecceities, or the Lack Thereof}
We now want to take a different tack and
examine some other motivations for considering non-Aristotelian logics, apart
from a desire to avoid the Godzilla problem.  We began by noting
that it is not good to build metaphysics into one's logic any more
than one should build facts of geography into trigonometry. On the
other hand, it is also desirable to have logics that are adequate to
the uses to which people frequently put language, and  to
the the kind of natural world we seem to find ourselves in.  We'll
point out that non-Aristotelian logics of some sort (not necessarily classical)
could well have applications in the logic of fiction, in philosophy,
and in physics itself.

\subsection{In Fiction}
We will not venture deeply into the logic and metaphysics of
fiction, except to note that it is not at all clear that simplistic
self-identity is deemed to hold for all entities treated in fiction,
the movies, and literature.  Walt Whitman, in his poem `Song of
Myself' \cite[p.\ 96]{Whitman}, famously stated, ``I contradict
myself\dots I contain multitudes.''  And what of a character caught
in a self-contradictory time loop in a bad science fiction story who
succeeds in committing suicide by shooting his grandfather?  At some
points along his worldline he exists if and only if he does not
exist. Only fiction, of course; but the point is simply that it is
open to authors even to question the self-identity of their
characters.

\subsection{In Philosophy}
Within the history of philosophy (and overlapping importantly into
physics) there has been and continues to be a debate between two
camps who have very different views about the metaphysics of time
and change. The Parmenideans see the world as static, the
Heracliteans see the world as inherently dynamic. Plato (to whose
work all of Western philosophy consists of footnotes, according to
Whitehead) proposed a synthesis of Parmenidean and Heraclitean
views: he distinguished between the world of Becoming (the unstable
physical or natural world) and the world of Being (the world of
stable ideal objects grasped by the intellect). Plato stated that
everything in the natural world, not only obviously changeable
things such as fire and water but even more apparently permanent
solid matter, was in a process of perpetual flux:
\begin{quote}
  Whenever we see anything in process of change, for example fire, we should
  speak of it not as \emph{being a thing} but as \emph{having a
  quality}\dots\ And in general we should never speak as if any of the things
  we suppose we can indicate by pointing and using the expressions `this
  thing' or `that thing' have any permanent reality:  for they have no
  stability and elude the designation `this' or `that' or any other that
  suggests permanence. \cite[p.\ 68]{Plato_Tim}
\end{quote}
Nothing in the world of Becoming is ever exactly a such-and-such and
thus one can never hope to fully grasp what it is.  Plato's words
are open to interpretation, of course, but his view seems to suggest
that objects in the natural world do not have sharp self-identity
because they are always in the process of becoming something else.

Nietzsche wrote in a similar vein:
\begin{quote}
Logic too depends on presuppositions with which nothing in the real world corresponds, for example on the presupposition that there are identical things, that the same thing is identical at different points of time\dots  \cite[\S 11]{Niet78}.
\end{quote}
On this view, even the claim that I am now sitting in the same chair that I sat in earlier this evening is (from the logical point of view) pure stipulation.  The later chair is like enough to the chair I sat on earlier for all practical purposes, so I might as well call it the same chair---but this is purely a stipulation justifiable only by its practical utility.

\subsection{In Physics}
We are certainly not suggesting that something is so just because
Plato or Nietzsche said it.  But such philosophical views, although imprecise and open to
interpretation, are not merely quaint relics of pre-scientific or 19th century romantic
thought; current professional debates on the reality of time and
change turn on the same Heralitean/Parmenidean point of dispute.

The static, plenum, or block-universe view is probably the
sentimental favourite of many recent physicists and philosophers of
science.  (Kurt G{\"o}del  was a  notable block universe theorist \cite{Yourgrau05}.)  However, because of the Indeterminacy
Relations and the fundamental non-Booleanity of quantum mechanics
(for an explanation of which see Bub \cite{Bub97}), it is not clear
that the static view is consistent with quantum mechanics
\cite{Peacock06}. The distinguished theorist Lee Smolin, trying to
understand the basis of the conceptual roadblocks which he insists
dog modern theoretical physics, remarks,
\begin{quote}
  I believe there is something basic we are all missing, some wrong
  assumption we are all making\dots\  I strongly suspect that the key
  is time.  More and more, I have the feeling that quantum theory
  and general relativity are both deeply wrong about the nature of
  time\dots\  We have to find a way to \emph{unfreeze} time---to
  represent time without turning it into space.  \cite[p.\ 257--7]{Smolin06}
\end{quote}
There is still no generally agreed upon method of unfreezing time,
but the questions of how to represent time, and whether or not time
is real or merely a funny sort of spatial dimension, are central
themes in current work in quantum gravity (the attempt to find a
quantum theory of spacetime structure).  If the \emph{logic} of
identity is to be of any use in talking about identity of objects in
time and space, it needs to be flexible enough to accommodate our
rapidly evolving picture of identity in the physical world.  Conceivably, some sort of non-Aristotelian logic as we conceive of it here could be a useful tool in Smolin's project to unfreeze physics.  

Paul Teller \cite{Teller98} has made some very relevant observations
about the way that \emph{haecceity}---the \emph{suchness} or
\emph{thisness} of an entity, that which presumably founds its
identity---is affected by quantum mechanics.  His explanation of
haecceity is very helpful:
\begin{quote}
  Traditionally, philosophy has talked about an object's
  ``haecceity'' to mark the idea that an object is distinct from
  all others in some manner that transcends all properties in any
  usual sense of the word `property.'  \dots\  let us take for
  granted some things that presuppose the applicability of strict
  identity:  that names can refer ``directly,'' that is without
  operating as definite descriptions; \emph{that repeated use of the same name
  picks out the same referent} [our emphasis]; that repeated use of
  the same variable bound by the same quantifier picks out the same
  referent; and that sets are defined extensionally\dots\ Now, a
  metaphysician might ask:  in virtue of what does strict identity
  apply to an object?  Haecceities\dots are supposed to be some
  metaphysical feature, principle, characteristic, or
  ``non-qualitative property'' which answers this question.
  \cite[p.\ 117]{Teller98}
\end{quote}
As Teller goes on to explain, in quantum statistics particles do
not have identities that can be tracked.  To adapt Teller's example,
if Bloggs has \$1000 in his bank account, it does not make sense to
ask, \emph{which} monetary  tokens (such as pennies) does he have \$1000 worth of?  All that
matters is that he is good for \$1000.  Similarly, if there are six
photons in a box, all this means is that we are good for six
photons;  it does not make sense to ask, \emph{which six photons
are in the box}?  As Teller explains, it has been found that if one
assumes that the photons have distinct, trackable identities the
way pennies do, one will count them wrong (because quantum particles
are permutation-invariant, unlike classical objects) and get the
wrong statistical predictions.  It is therefore highly questionable
that quantum mechanics allows for the notion of haecceity (and
thereby self-identity) in anything like the classical sense.  (See \cite{FK2006}  for a detailed exploration of this problem.)

There is another respect in which quantum mechanics suggests
something like the Heraclitean view. Any assertion in quantum
physics has to be operationally grounded; by what measurement
procedure could we know that a particle is identical to itself?
Well, we might have to interact with it twice, and quantum mechanics
tells us that there is no clear meaning to saying that if we measure
(say) an electron at a certain spot, and then a tiny fraction of a
second later measure another electron at nearly the same spot, that
we have detected the same electron that we detected in the first
measurement.  To adapt a famous phrase from Heraclitus, we do and do
not observe the same electron.

At around this point in the discussion, classically-minded thinkers are sometimes moved to exclaim, ``But
dammit, everything just \emph{is} identical to itself!'' This is an example of what can be
called the \emph{table-pounding argument} because such statements are
often accompanied by pounding on a convenient mid-sized object.
Unfortunately for the classical realist, quantum mechanics remains
unmoved by any amount of furniture-thumping. Modern physics
certainly suggests, and arguably demands, that we
live in an extreme Heraclitean world of flux where
self-identity cannot be asserted, or at least cannot be asserted in
a classical way---except as a convenient approximation at scales where quantum effects can be ignored.  

One does not necessarily employ first order predicate logic to
reason about quantum mechanics, but in order for logic to be as useful
as possible, in the kind of physical world we live in,  it should be
equipped to express facts of quantum mechanics as required and
should therefore not have in-built assumptions that would conflict
with quantum mechanics. In the spirit of Putnam's recommendations
\cite{Putnam68}, it is desirable to seek a way of doing classical
logic that would naturally generalize to quantum logic.

To conclude this section, we quote a favourite story from Bertrand Russell:
\begin{quote}
  It is obvious that, if you think of all the things that are in the world,
  they cannot be divided into two classes---namely,  those that exist, and
  those that do not.  Non-existence is, in fact, a very rare property.
  Everybody knows the story of the two Germanic pessimistic philosophers, one
  of whom exclaimed:  `How much happier were it never to have been born.' To
  which the other replied with a sigh: `True!  But how few are those who
  achieve this happy lot.'  \cite[p. 147]{Russell56}
\end{quote}
We have a great deal of respect for Russell, but it is by no means
obvious that existence is such a clear-cut
concept in a quantum universe. Perhaps if Teller is right then non-existence is not such a
difficult property to attain or at least to approximate after all,
at least within the limits allowed by the Uncertainty Principle.

To summarize:  self-identity holds for some idealized objects, such
as the natural numbers, and it is approximately enough true to not
be misleading for many mid-sized physical objects such as tables and
chairs.  Physics tends to suggest that it could well be simply dead
wrong at scales where quantum mechanics is important, although this
remains an important open question.  But again, if predicate logic
is to be as widely applicable as possible to reasoning about things
in the nature world, not to mention fictional objects, it must not
be burdened with the \emph{presumption} that everything that can be
quantified over is necessarily self-identical.

\section{The Null Object}
Now we will show how non-Aristotelian even a nearly-classical open logic has to be if it is to have enough expressive power to be useful in mathematics and daily reasoning.

Let us consider the following variation on the Categorical Godzilla proof, which we dub the \emph{Conditional Godzilla}:
\begin{tabbing}
mm \= mmmm \= mmm  \= mmmmmmmmmmmmmm  \= mmmmmm   \kill
    \>1     \>(1)  \> $\forall x(x=x)$           \>A  \\
    \>1     \>(2)  \> $\text{Godzilla} = \text{Godzilla}$      \>1 UE  \\
    \>1     \>(3)  \> $\exists x(x = \text{Godzilla})$    \>2 EI
\end{tabbing}
Thus by conditional proof we have
\begin{equation}
 \ts \forall x(x=x) \onlyif \exists x(x = \text{Godzilla}).
\end{equation}
Here we have assumed $\forall x(x=x)$ rather than taken it as a
theorem. We're on firmer ground in that respect.  But we still end
up with a peculiar result: assuming the self-identity of all
objects in the universe of discourse apparently also allows us to
prove the existence of anything in that universe, only this time not as an ersatz theorem but as a consequence of the assumption on line (1).   So if we give all entities in $\mathcal{U}$ the benefit of the doubt and grant them self-identity, we are still committed to their existence---but now the dependence of the existence result upon assumption is obvious. 

In part this proof is simply an illustration of the point noted earlier, that in first order predicate logic, to posit any property (including
self-identity) of an entity is to imply that the entity exists. So if we give all entities in $\mathcal{U}$ the benefit of the doubt and grant them self-identity---not as a presumed logical or metaphysical truth but simply for the sake of argument---that still implies that they exist.  On the other hand, non-existent things are
non-self-identical, simply because no properties of any sort can \emph{in fact}  be predicated
of them at all (regardless of how they were defined).  But what happens if the facts of a matter demand that we deny the existence of something that we had provisionally admitted to $\mathcal{U}$?

Continue the above proof as follows:
\begin{tabbing}
mm \= mmmm \= mmm  \= mmmmmmmmmmmmmmmmmmm  \= mmmmmm   \kill
    \>     \>(4)  \> $\forall x(x=x) \onlyif \exists x(x = \text{Godzilla})$           \>1--3 CP  \\
    \>5     \>(5)  \> $-\exists x(x = \text{Godzilla})$           \>A \text{(An empirical given.)}  \\
    \>5     \>(6)  \> $-\forall x(x=x)$           \>4,5 MT  \\
    \>5     \>(7)  \> $\exists x(x \neq x)$           \>6 Duality  
\end{tabbing}
By asserting on line (5) the empirical fact that a certain described entity does not exist, we seem to be forced into a bizarre existence claim anyway!  One approach could be to simply not make assertions like (5) on the grounds that in classical predicate logic we hold or pretend that all names refer.  But predicate logic would be greatly restricted in its usefulness if it could not assert the non-existence of a putative entity known only by a name or description.  That is something that we are entirely free to do in ordinary language, as well as scientific and mathematical reasoning.  As we have already noted, one of the most useful tools of reasoning is the ability to discuss something hypothetically, be it the largest prime or the gunman on the grassy knoll.  What we need is a natural interpretation of the odd thing whose existence is cited in line (7).  We suggest that it may be useful to think of this object, or ``object'', as a \emph{null entity}.   In a logical system it acts in a way analogous to the ground in an electrical circuit; it is  the elephant graveyard for all names and descriptions which fail to refer.  

One can see a foreshadowing of this approach in Carnap's interpretation of Frege's solution to the problem of improper definite descriptions \cite[35--9]{Carnap47}.  As Carnap explains, Frege was concerned to construct his ideal logical language so that a definite description picked out a unique object. There is an obstable in the cases of improper definite descriptions, terms which have the form of a description but which name nothing or many things.  Carnap observes that a possible response is ``to count among the things also the \emph{null thing}, which corresponds to the null class of space-time points'' \cite[36]{Carnap47}.   It is beyond the aim of this paper to go into more detail regarding the problem of improper definite descriptions, but it is well worth noting that there is a precedent in the literature for this kind of solution to problems not too dissimilar to those in which we are most interested here.

For Carnap, a null object $a_0$ is simply a name that is left free---it does not denote anything.  Since names can denote anything we want, it is open to us to simply leave one name unassigned in the course of a piece of reasoning.  Classically (i.e., in a logic with a sharp concept of identity), a null object can be naturally defined as follows:
\begin{equation}
a_0 :=  x(x \neq x).  \label{NullDf}
\end{equation}
(``Let $a_0$ be an $x$ such that $x \neq x$.'')\footnote{
Our notation here is non-standard and requires some explanation.  It would be more common to write a definition such as this using Hilbert's $\varepsilon$-symbol, as $a_0 :=  \varepsilon x(x \neq x)$.  However, the usual reading of $\varepsilon xA$ commits us to more than we want when discussing null objects:  Leisenring \cite[p.\ 1]{Leis69} says, ``Intuitively, the $\varepsilon$-term $\varepsilon xA$ says `an $x$ such that if anything has the property $A$, then $x$ has that property'.''  Our simpler notation is inspired by the reading of $\exists xFx$ as ``There exists an $F$.''  The symbol $\exists$ is ``there exists,'' and $xFx$ is ``an $F$''.  Any expression of the latter form can be called an \emph{indefinite denotator}.  Precisely because it is so thoroughly indefinite, it could (depending on the facts of the matter) denote nothing.  
}
It is not essential that $a_0$ be defined this way:  that is, as whatever is non-self-identical.  The notion of a null object could easily survive a loosening or broadening of the concept of self-identity.  However, if we accept the classical notion of the null object, then lines (4)--(7) above do not in fact demonstrate the existence of anything at all---precisely because they demonstrate the existence of the null object, which is not anything at all.  So we have not stumbled into an existence claim simply because we tried to deny the existence of something.  

The move to free logic provides another natural motivation for considering the null object.  Again, a free logic is defined as a system of predicate logic which allows
for the possibility of empty domains of discourse and predicates
that do not refer.  Consider $Fa \vel -Fa$.  This is true even in
$\emptyset$ (the null set, an empty universe) because if a term $a$
does not refer then $-Fa$ holds for any predicate $F$.  (If
the present King of France doesn't exist then it is true that he is not
bald.)  However, the fact that $Fa \vel -Fa$ is true in an empty
universe seems to license a dubious inference:
\begin{tabbing}
mm \= mmmm \= mmm  \= mmmmmmmmmmmmmmmmmmm  \= mmmmmm   \kill
    \>     \>(1)  \> $Fa \vel -Fa$           \>Theorem  \\
    \>     \>(2)  \> $\exists x(Fx \vel -Fx)$           \>1 EI (???)
\end{tabbing}
That is, we seem to have once again inferred a theorem asserting the existence of an entity---this time, an apparent element of $\emptyset$!

One way to deal with this is to not talk about empty universes; and
this is what is usually done in elementary predicate logic texts,
where the puzzle of the empty universe is either not mentioned or
glossed over. Another way to block this inference is to
not allow EI for empty universes, but this requires that one know in
advance that a universe is empty, and we are supposed to be able to
do logic without any existential presumptions at all.   Here, again, a null object can help us.  Classically, we have $\emptyset = \{x|x \neq x\}$.  Then $a_0 \in \emptyset$.   (Indeed, it's the only element of $\emptyset$.)  Then line (2) above can only apply to $a_0$ (in $\emptyset$), and we have (in $\emptyset$) $Fa_0 \vel -Fa_0$.  Now, $Fa_0$ is false for any predicate $F$; by bivalence $-Fa_0$ and therefore $Fa_0 \vel -Fa_0$ are true in $\emptyset$.  So the above deduction is valid in $\emptyset$.  

Please note:  this bit of reasoning does not imply that there actually is anything at all \emph{in} $\emptyset$!

Null objects thus allow us to extend the validity of certain puzzling arguments to empty universes in a natural way, and allows us to express the non-existence of named or described entities when the facts demand that we do so---so long as we decide that we can live with one
more very odd sort of mathematical creature under the floorboards of
everyday reasoning.  Just as set theory would be hobbled without the
formal device $\emptyset$, and arithmetic could not operate without 0, it could well be that predicate logic has been hobbled all
along  without a formal, placeholder referent for names and descriptions that do not refer.  

One further point:  is the null object self-identical?    Clearly not, by its very definition, Eq.\ (\ref{NullDf}).  So even if we want to keep our predicate logic as close to classical as possible (by adhering to =I as a global assumption but making no other changes in our deductive methods), if we also want to be able to assert that some names or descriptions fail to refer there must be at least one non-self-identical (and therefore null) object in every domain of reference.   To this extent, then, even an open classical logic is non-Aristotelian.  

We were tempted to speak of \emph{the} null object, but Gillman Payette (private communication) was quick to point out to us that even this would be saying too much about it.  Suppose that we tried to define $a_0$ as follows:
\begin{equation}
a_0 := \Russell x(x \neq x).  \label{TheNullDf}
\end{equation}
(``Let $a_0$ be \emph{the} $x$ such that $x \neq x$.'')  To speak of anything as  \emph{the} $F$ is to allow that it can be equal to something bearing a proper name.  Suppose that $a_0 = b$.  So long as we are allowed =E (and we could not do much useful reasoning about identity without it), we can substitute $a_0$ for $b$ and get that $a_0$ is self-identical---precisely the thing we don't want.  So $a_0$ cannot even have the property of uniqueness.  In this respect the analogy between null objects, and the null set and 0, breaks down because the latter entities can stand in identity relations.  So the null object or objects must remain utterly indefinite; our $a_0$ is just a placeholder for the absence of all properties, demanded by the syntax of
predicate logic.  

There is one further intriguing observation to be made about null objects.  The definition (\ref{NullDf}) is a very natural way to specify a null object in Aristotelian logic.  But suppose we want to consider non-Aristotelian logics where classical identity is not always available.  We would need a more general conception of null objects.  If we are allowed to quantify over predicates, then we could define
\begin{equation}
  a_0 := x(\forall F(-Fx)).  \label{NullDfSteroids}
\end{equation}
This has the advantage that it could apply to logics without classical identity.  But one encounters a challenge that should by now be very familiar to logicians.  Let us say that an object is \emph{prediphobic} if it will admit of no predicates whatsoever:
\begin{align}
  Px &:= \text{$x$ is prediphobic} \\   \nonumber
       & := \forall F(-Fx)  \label{PrediphobicDf}
\end{align}
Then clearly $Pa \onlyif -Pa$, and given $Pa$ as well, we have detonation.  This is a typical instance of the hazards of second-order logic.  It also suggests that any logic in which classical identity fails, but in which a more general null object of the form (\ref{NullDfSteroids}) is desired, would have to be paraconsistent in some sense.

\section{Summing Up}
The notion that =I is a theorem or logical truth leads to the unacceptable result that the existence of any named or described entity can be proven as a matter of logic.   The most natural way out of this embarrassment is to think of =I as an assumption, not a ``logical truth''.  Indeed, there are ample reasons within fiction, philosophy, and physics why we might want to speak of entities whose self-identity is in doubt, and we should have logics that are open to this possibility.  We suggest that a logic in which =I holds for all non-null entities in
its domain be called \emph{Aristotelian}; otherwise, \emph{non-Aristotelian}.  If =I is taken to hold as a logical truth, we'll call that a classical Aristotelian logic.  (We expect that such logics will sooner or later become historical curiosities.)  If =I is taken to hold merely as an assumption, but still an assumption applying to all objects in the domain of discourse, we'll say that such a logic is an \emph{open} Aristotelian logic.  We conjecture---though this remains to be shown in full rigour---that an open Aristotelian logic can do all the deductive work that classical Aristotelian logic can do, without falling into the absurdity of the Categorical Godzilla.  Beyond this, an important research project is to explore possible non-Aristotelian logics.  

Our attempt to de-ontologize logic (by removing =I as a theorem) ironically forces us to include null objects in any domain of reference, even when the logic is Aristotelian.  But this is not an expansion of our ontology (except for a modest addition to our collection of symbols) because the null object is not any \emph{thing}, even though it may turn out to be just as indispensable as $\emptyset$ and 0.

One result seems clear:  we can no more take it to be a logical truth that everything is self-identical than we can take it to be a logical truth that everything is green.  If we do impute self-identity to all non-null members of a domain of discourse, it is only by courtesy or because it is a domain (such as $\mathbb{N}$ or the furniture in someone's office) for which we have good reason to think that self-identity holds throughout.  And any logic that hopes to be adequate to the ordinary demands of discourse in the real world must always allow for the possibility that self-identity fails for some entities of interest.

\section*{APPENDIX:  Some Properties of an Open Classical Logic}
We have a lot of work to do in order to clarify the properties of
open classical logic, let alone explore other non-Aristotelian logics that might
be feasible and link them with free logics. Here we list without
proof some immediate consequences and properties of an open classical logic.

It is simply a result of classical first order logic that we have $a=a\rightarrow\exists x(x=a)$. For a natural deduction system, this is very easily provable by the rules of EI and conditional introduction, with $a=a$ as an assumption. If the system is stated in a normal axiomatic presentation, we have as an axiom that $B\rightarrow\exists xB$, where some instances of a particular name $a$ occurring in $B$ can be replaced by the variable $x$, bound by the existential quantifier. The identity $a=a$ is such a formula, and so $a=a\rightarrow\exists x(x=a)$ is just an instance of this axiom. A result of this with a classical theory of identity is that $\exists x(x=b)$ is provable for any name $b$, whereas in our open logic, $\exists x(x=b)$ follows only in cases where we explicitly assume $b=b$. 

We retain the provability of sequents of the form $a=b,Pa\vdash Pb$. In the natural deduction setting, this is enforced by a rule of identity elimination. More generally, it is a form of Leibniz's law. For axiomatic and sequent calculus purposes, we can simply include this sequent as a primitive rule, as is common in the proof theory literature. 

In a sequent system our proposal amounts to rejecting $\vdash a=a$ as an axiom. However, given the inclusion of the identity elimination rule, all one need to do to reason more or less usually with identity is to include a premise of the form $a=a$; that is to include $a=a$ on the left-side of the $\vdash$. Including this extra premise is always admissible, by the rule of thinning, and so, just as in the natural deduction case, we demand only that one make one's \emph{extra-logical} assumptions about identity explicit.

A result of this rule is that $a=b\vdash a=a$ is provable. So, under the assumption that $a$ is equal to anything, we have that $a=a$, and if $a\neq a$, then $a$ is identical to nothing. With identity elimination, we can also easily prove that $a=b, b=c\vdash a=c$, where $b=c$ is the formula in which $a$ is substituted for $b$ to attain the conclusion. Hence, we can clearly prove that $(a=b\land b=c)\rightarrow a=c$. Similarly, we can show that $a=b\equiv b=a$. So, as a result, $=$ enjoys symmetry and transitivity in all cases, and reflexivity in those cases where it applies at all. While this `conditional' reflexivity is strictly weaker than reflexivity, it guarantees that in contexts where we have assumed self-identity to hold of the names we reason with, identity behaves classically. The differences are, of course, with those names of which we make no such assumption.

Of course, since we reject that $a=a$ is a theorem, we do not have that the formula $a\neq a$ allows for the proof of any formula whatsoever. In general, the assumption $a\neq a$ will only generate triviality when we explicitly assume some other formula which implies $a=a$, because we retain the rule of \emph{explosion} (ex falso quodlibet).  

We also retain the \emph{law of excluded middle}, because the propositional fragment of our logic is purely classical. So, we have that $\vdash a=a\lor a\neq a$; however, as we do not have that $\vdash a=a$ or $\vdash a\neq a$, the logic is not \emph{prime}. This is just to say that it is not universally true that $\vdash A\lor B$ holds only when either $\vdash A$ or $\vdash B$ holds. This leaves a potentially interesting avenue to intuitionist open logic available. It strays beyond our aims to investigate such a logic, but we note here that such an approach may have interesting consequences for common subjects which motivate intuitionistic logic, for instance mathematics\footnote{Consider the fact that from ZFC alone we can prove neither $\aleph_1=2^{\aleph_0}$ nor $\aleph_1\neq 2^{\aleph_0}$. Of course, giving a set theoretic analysis of $=$ which matches our assumptions about this predicate goes beyond our purposes, but it may be a valuable way forward. This is in contrast to the well-understood notion of \emph{bijection} underwriting $=$ in ZFC, and other common set theories.} and other areas where epistemic restrictions on our knowledge are salient.

No proposal to amend an established logic can be taken seriously
until the metatheory of the amended logic is worked out.  We have
yet to do this. It seems very likely that \emph{dropping} a rule of
inference leaves us on safe ground with regard to soundness.  In
particular, dropping =I as a theorem makes no difference at all to
what can actually be deduced with first order predicate logic with
identity except that we lose certain inferences that (for reasons we
have explained) we would like to lose anyway.

Completeness is a more difficult question: to show that open
classical first order logic with identity is complete we would have
to show that any formula we can \emph{no longer} prove (by having
removed =I as a rule of deduction) is false in some models that
leave true all the theorems that did not require =I.  To put it
another way, completeness is all about whether one can prove all of
the tautologies in a system with the resources of the system.  The
Categorical Godzilla shows us that if the standard approach is
complete and sound, $\exists x(x=a)$ must be a tautology for any
term $a$ in every possible model; and as we have noted this could
make sense only if it is somehow known that every term in the
language refers. In an open classical logic $\exists x(x=a)$ is most certainly
not a tautology in general, but rather a statement that depends upon
the facts of a case.  So again, we think that dropping =I as a rule
of inference will only prevent one from being able to prove formulas
that have no business being tautologies anyway.  But this question
requires a more thorough study.

If nothing else, this appendix makes it clear that this work is a very first step into an interesting new territory, about which almost everything is as yet unknown. Our mere hope at this stage is that the reader is intrigued by our proposal, and convinced enough by our philosophical argumentation to think that it may be worth developing in more detail.

\section*{Acknowledgements}
We are grateful to the Universities of Lethbridge and Connecticut for material and financial support, and to Bryson Brown for many valuable discussions about logic and other matters.  Karl Laderoute pointed us to the remarks by Nietzsche, and John Woods has provided helpful comments about the notion of identity.  Gillman Payette commented insightfully on an earlier version at the Canadian Philosophical Association at Congress 2016 (at the University of Calgary).  We alone are responsible for any errors or omissions that may remain in this paper despite all the good advice we've received.

\newpage
{\frenchspacing\raggedright\small\linespread{1.25}\selectfont\raggedright

}

\end{document}